\documentclass[11pt]{article}
\usepackage{tipa}
\usepackage{stmaryrd}

\usepackage{amsmath, amssymb, latexsym}
\usepackage{psfrag,epsfig,amsfonts,amsmath,latexsym,amsthm,amssymb,amscd }

\usepackage{color,amsmath,amssymb, amsfonts,amstext,amsthm}

\usepackage{epsfig, graphicx, graphics}
\usepackage{bbm}
\usepackage[colorlinks]{hyperref}
\usepackage{mathrsfs}
\usepackage{url}

\newcommand{\eps}{\varepsilon}

\newcommand{\EX}{{\mathbb{E}}}
\newcommand{\PX}{{\mathbb{P}}}

\newcommand{\R}{\mathbb{R}}

\theoremstyle{remark}
\newtheorem{remark}{Remark}
\newtheorem{example}[remark]{Example}
\topsep=\parskip

\title{A parameter estimation method based on  random slow manifolds
 \footnote{This work was done while Jian Ren was visiting the Institute for Pure and Applied Mathematics (IPAM), Los Angeles,   USA.
It was partly supported by the NSF Grant  1025422 and   the NSFC grant 11271290. } }

\author{Jian Ren$^1$ and Jinqiao Duan$^2$
\\1. School of Mathematics and Statistics\\
Huazhong University of Science and Technology\\ Wuhan, 430074, China. Email: renjian0371@gmail.com \\
\&
\\2. Institute for Pure and Applied Mathematics (IPAM), University of California\\Los Angeles, CA 90095, and \\
Department of Applied Mathematics, Illinois Institute of Technology\\ Chicago, IL 60616,  USA. Email: duan@iit.edu   }

\date{\today}

\begin{document}
\maketitle

\begin{abstract}
A parameter estimation method is devised  for a
slow-fast stochastic dynamical system, where often only the slow component is observable. By using the observations  only on the slow
  component,  the system parameters are estimated by working on the slow system on the random slow manifold. This offers a benefit of dimension reduction in quantifying parameters in stochastic dynamical systems. An example is presented to illustrate this method, and verify that the parameter estimator based on the lower dimensional, reduced slow system is a good approximation of the parameter estimator for  original slow-fast stochastic dynamical system.


\end{abstract}




\section{Introduction}

Invariant manifolds provide geometric structures for understanding
dynamical behavior of nonlinear systems under uncertainty. Some
systems   evolve on fast and slow time scales, and may   be
modeled by  coupled singularly perturbed stochastic ordinary   differential equations (SDEs). A slow-fast stochastic system may have a special invariant manifold called a random slow manifold that capture the slow dynamics.

 We consider a stochastic slow-fast system
\begin{eqnarray}\label{Slow-Equation-Stoch000}
&&\dot{x}=A{x}+  f(x, y),
\quad  x(0)=x_0 \in \mathbb{R}^n,\\
\label{Fast-Equation-Stoch000}
&&\dot{y}=\frac{1}{\eps}B{y}+\frac{1}{\eps}g(x, y)+
\frac{\sigma}{\sqrt{\eps}}\dot{W}_t, \quad y(0)=y_0\in\mathbb{R}^m,
\end{eqnarray}
where $A$ and $B$ are  matrices, $\eps$ is a small positive
parameter measuring slow  and fast scale separation,
  $f$ and $g$ are nonlinear Lipschitz continuous functions with Lipschitz constant $L_f$ and $L_g$ respectively, $\sigma$ is a noise intensity constant, and $\{W_t:t\in\mathbb{R}\}$
is a two-sided $\mathbb{R}^m$-valued Wiener process (i.e., Brownian motion) on a probability space $(\Omega,  \mathcal{F}, \PX)$.
   Under a gap
condition and for $\eps$ sufficiently small, there exists a random
slow manifold $y=h^\eps(\xi, \omega)$, $\omega \in \Omega$, as in \cite{Schm, FuLiuDuan},
for slow-fast stochastic system
\eqref{Slow-Equation-Stoch000}-\eqref{Fast-Equation-Stoch000}. When
the nonlinearities $f, g$ are only locally Lipschitz continuous but
the system has a random absorbing set (e.g., in mean-square norm),
we conduct a cut-off of the original system. The new system will
have a random slow manifold which captures the original system's
slow dynamics.

 The random slow manifold is the graph of a random nonlinear mapping
   $h^\eps(\xi, \omega)=\sigma\eta^\varepsilon(\omega)+~\tilde h^\varepsilon(\xi, \omega)$,  with $\tilde h^\varepsilon(\xi, \omega)$ determined by a
Lyapunov-Perron integral equation \cite{FuLiuDuan},
$$\tilde h^\eps(\xi, \omega)=\frac1{\eps}\int_{-\infty}^0 e^{-\frac{B}{\eps}s}g( x(s, \omega, \xi), y(s, \omega, \xi)+\sigma \eta^\eps(\theta_s\omega))\,ds,$$
with $\eta^\eps(\theta_t\omega)=\frac1{\sqrt{\eps}}\int_{-\infty}^t
e^{\frac{B}{\eps}(t-s)}\,dW_s$ and
$\eta^\eps(\omega)=\frac1{\sqrt{\eps}}\int_{-\infty}^0
e^{-\frac{B}{\eps}s}\,dW_s$.
The random slow manifold
exponentially attracts other solution orbits. We will  find an analytically approximated random slow manifold for
sufficiently small $\eps$, in terms of an asymptotic expansion in
$\eps$, as in \cite{SunDuan1, SunDuan2}. This slow manifold may also be numerically computed as in \cite{Kan2013}.  By restricting to the slow manifold, we obtain a lower
dimensional reduced system of the original slow-fast system
\eqref{Slow-Equation-Stoch000}-\eqref{Fast-Equation-Stoch000}, for
$\eps$ sufficiently small
  \begin{eqnarray} \label{slowdynamics}
  \dot{x}=A{x}+ f(x, \tilde h^\eps(x, \theta_t\omega)+\sigma \eta(\psi_\varepsilon\omega)),\quad x\in
\;\;\mathbb{R}^n,
  \end{eqnarray}
where $\theta_t$ and $\psi_\varepsilon$ are defined in the next section.


If the original slow-fast system
\eqref{Slow-Equation-Stoch000}-\eqref{Fast-Equation-Stoch000} contains unknown system parameters, but only the slow component $x$ is observable, we conduct parameter estimation using the slow system
\eqref{slowdynamics}. Since the slow system is lower dimensional than the original system, this method offers an advantage in computational cost, in addition to the benefit of using only observations on slow variables.

This paper is arranged as follows. In the next section, we obtain an approximated random slow manifold and thus the random slow system. Then in Section 3, we present a method for parameter estimation on the slow manifold. Finally, a simple example is presented in Section 4 to illustrate our method.

\section{Random slow manifold and its approximation}\label{slow manifold-approximation}

By  a random transformation
\begin{eqnarray}\label{random transformation}
\begin{matrix}\begin{pmatrix}  X \\
Y
\\\end{pmatrix} :=\mathcal V_\varepsilon(\omega,x, y)=
\begin{pmatrix}  x \\y-\sigma\eta^\varepsilon(\omega)
\\ \end{pmatrix} \end{matrix},
\end{eqnarray}
we convert the SDE
system \eqref{Slow-Equation-Stoch000}-\eqref{Fast-Equation-Stoch000}
to the following   system with random coefficients
\begin{eqnarray}\label{Slow-Equation-random000}
&&\dot X(t)=AX(t)+ f(X(t), Y(t)+ \sigma
\eta^\eps(\theta_t\omega)),\\
\label{Fast-Equation-random000} &&\dot
Y(t)=\frac{1}{\eps}BY(t)+\frac{1}{\eps}g( X(t), Y(t) +\sigma
\eta^\eps(\theta_t\omega)),
\end{eqnarray}
where $\eta^\eps(\omega)=\frac1{\sqrt{\eps}}\int_{-\infty}^0
e^{-\frac{B s}{\eps}}\,dW_s$ is the stationary solution of linear
system $dy^\eps = \frac{B}{\varepsilon}y^\eps dt +
   \frac{\sigma}{\sqrt{\varepsilon}}dW_t$.
 Here $ \theta_t: \Omega \to \Omega $ is the Wiener shift implicitly defined by
 $W_s(\theta_t\omega)= W_{t+s}(\omega)-W_t(\omega)$.
   Note that
$\eta^\eps(\theta_t\omega)=\frac1{\sqrt{\eps}}\int_{-\infty}^t
e^{\frac{B}{\eps}(t-s)}\,dW_s$.

Define a mapping (between random samples) $\psi_\varepsilon:\Omega\rightarrow\Omega$ implicitly by
$~W_t(\psi_\varepsilon\omega) =
\frac{1}{\sqrt{\varepsilon}}W_{t\varepsilon}(\omega)$. Then
$\frac{1}{\sqrt{\varepsilon}}W_{t\varepsilon}(\omega)$ is also a
Wiener process with the same distribution as $W_t(\omega)$.
Moreover, $\eta^\varepsilon(\theta_{t\varepsilon}\omega)$ and
$\eta^\varepsilon(\omega)$ are identically distributed with
$\eta(\theta_t\psi_\varepsilon\omega)=\int_{-\infty}^t
e^{B(t-s)}\,dW_s(\psi_\eps \omega)$ and
$\eta(\psi_\varepsilon\omega)=\int_{-\infty}^0 e^{-B
s}\,dW_s(\psi_\eps \omega)$, respectively.

By a time change $\tau=t/\varepsilon$ and using the fact that
$\eta^\eps(\theta_{\tau\varepsilon}\omega)$ and
$\eta(\theta_\tau\psi_\varepsilon\omega)$ are identically distributed, the system
\eqref{Slow-Equation-random000}-\eqref{Fast-Equation-random000}  is reformulated as
\begin{eqnarray}\label{Slow-Equation-random000fast}
&&X'= \varepsilon [AX+ f(X,
Y+ \sigma \eta(\theta_\tau\psi_\varepsilon\omega))],\\
\label{Fast-Equation-random000fast}&& Y'=BY+g( X, Y +\sigma
\eta(\theta_\tau\psi_\varepsilon\omega)),
\end{eqnarray}
where $'=\frac{d}{d\tau}$.


 We make the following two hypotheses.

 $\mathbf{H1}$: There are positive constants
$\alpha$, $\beta$  and $K$,  such that for every $x\in \mathbb{R}^n$
and $y\in \mathbb{R}^m$,  the following exponential estimates hold:
 $$|e^{At}x|_{\mathbb{R}^n} \leq K e^{\alpha
t}|x|_{\mathbb{R}^n},\quad t\leq 0;\quad \quad
|e^{Bt}y|_{\mathbb{R}^m} \leq K e^{-\beta t}|y|_{\mathbb{R}^m},\quad
t \geq 0.$$

$\mathbf{H2}$: $\beta>K L_g$.\\

Then there exists  a random slow manifold $\mathcal{\tilde M}^\eps(\omega)=\{(\xi,
\tilde h^\eps(\xi, \omega)):\xi\in \mathbb{R}^n\}$ for the random
system
\eqref{Slow-Equation-random000}-\eqref{Fast-Equation-random000},
with $\tilde h^\eps$ being
expressed as  follows   \cite{FuLiuDuan, RenDuanJones},
$$\tilde h^\eps(\xi, \omega)=\frac1{\eps}\int_{-\infty}^0 e^{-\frac{B}{\eps}s}
g(X(s, \omega, \xi), Y(s, \omega, \xi)+\sigma
\eta^\eps(\theta_s\omega))\,ds.$$ We  can get a small $\eps$
approximation for $\tilde h^\eps$.  Start with   the expansion
$Y=Y_0+\varepsilon Y_1+ \mathcal{O}(\varepsilon^2)$ and the integral
expression   $X=X(0) + \varepsilon\int_0^\tau [AX+ f(X, Y+ \sigma
\eta(\theta_r\psi_\varepsilon\omega))]\,dr $  in the   system
\eqref{Slow-Equation-random000fast}-\eqref{Fast-Equation-random000fast}.
Using the Taylor expansions of $f(X, Y+ \sigma
\eta(\theta_\tau\psi_\varepsilon\omega))$ and $g(X, Y+ \sigma
\eta(\theta_\tau\psi_\varepsilon\omega))$ at point $(x_0, Y_0+\sigma
\eta(\theta_\tau\psi_\varepsilon\omega))$ and denoting $\xi=~X(0)$,
we obtain
\begin{eqnarray} \label{slowmanifold}
  \tilde h^\eps(\xi, \omega)=~\tilde
h_d^\eps(\xi, \omega)+~ \eps \tilde h_1^\eps(\xi, \omega) +~
\mathcal{O}(\eps^2),
\end{eqnarray}
 for $\eps $ sufficiently small, with
\begin{eqnarray}\label{h_d}
\tilde h_d^\eps(\xi, \omega)=\int_{-\infty}^0 e^{-Bs}g(\xi,
Y_0(s)+\sigma\eta(\theta_s\psi_\varepsilon\omega))\,ds,
\end{eqnarray}
and
\begin{eqnarray}\label{h_1}
\tilde h_1^\eps(\xi, \omega)&=&\nonumber\int_{-\infty}^0
e^{-Bs}\Big\{g_x( \xi,
Y_0(s)+\sigma\eta(\theta_s\psi_\varepsilon\omega))\big[As\xi+\int_0^s
 f(\xi,
Y_0(r)+\sigma\eta(\theta_r\psi_\varepsilon\omega))\,dr\big]\\{}&&\quad\quad\quad\quad\quad\quad
+g_y( \xi,
Y_0(s)+\sigma\eta(\theta_s\psi_\varepsilon\omega))Y_1(s)\Big\}\,ds,
\end{eqnarray}
where $Y_0(t)$ and $Y_1(t)$ satisfy the following random
differential equations, respectively
\begin{eqnarray}\label{y_0}
\begin{cases}
Y'_0(\tau)=BY_0(\tau)+g( \xi,
Y_0(\tau)+\sigma\eta(\theta_\tau\psi_\varepsilon\omega)), \\
Y_0(0)=h_d^\eps(0, \omega),
\end{cases}
\end{eqnarray}
and
\begin{eqnarray}\label{y_1}
\begin{cases}
Y'_1(\tau)=\big[B+g_y( \xi,
Y_0(\tau)+\sigma\eta(\theta_\tau\psi_\varepsilon\omega))\big]Y_1(\tau)\\
\quad\quad\quad\quad+g_x(\xi,
Y_0(\tau)+\sigma\eta(\theta_\tau\psi_\varepsilon\omega))\big\{A\tau\xi+\int_0^\tau
f(\xi, Y_0(s)+\sigma\eta(\theta_s\psi_\varepsilon\omega))\,ds\big\} , \\
Y_1(0)= h_1^\eps(0, \omega).
\end{cases}
\end{eqnarray}

Noticing that
 $\eta^\eps(\omega)$ and $\eta(\psi_\varepsilon\omega)$ have  identical distribution, together with the fact that  $$\eta^\eps(\omega)=\frac1{\sqrt{\eps}}\int_{-\infty}^0
e^{-\frac{B s}{\eps}}\,dW_s=\frac1{\sqrt{\eps}}\int_{-\infty}^t
e^{\frac{B (t-u)}{\eps}}\,dW_u=\eta^\eps(\theta_t\omega), \quad
u=s+t,$$
we see that $\eta^\eps(\theta_t\omega)$ and
$\eta(\psi_\varepsilon\omega)$ are identically distributed. Recall
  the transformation introduced in the beginning of
this section $X(t)=x(t)$ and
$Y(t)=y(t)-\sigma\eta^\varepsilon(\theta_t\omega)$. We thus obtain a
lower dimensional,   slow stochastic system on the slow manifold,
\begin{eqnarray}
\dot x &=& Ax + f(x, \tilde h^\eps(x, \theta_t\omega)+\sigma
\eta(\psi_\varepsilon\omega)),  \label{sloweqn}
\end{eqnarray}

Moreover,
 $
\hat h^\eps(\xi, \omega)= \tilde h_d^\eps(\xi,
\omega)+ \eps \tilde h_1^\eps(\xi, \omega)
 $
is   an approximation or a first order truncation of $\tilde h^\eps(\xi, \omega)$,
and hence we have an approximate slow    system
\begin{eqnarray}
\dot x &=& Ax + f(x, \hat h^\eps(x, \theta_t\omega)+\sigma
\eta(\psi_\varepsilon\omega)).  \label{sloweqn-app}
\end{eqnarray}
 This is the \emph{slow system} we will work on for parameter estimation in the next section.

\section{Parameter estimation on a random slow manifold}

If the original slow-fast system
\eqref{Slow-Equation-Stoch000}-\eqref{Fast-Equation-Stoch000} contains an unknown system parameter $a$, but only the slow component $x$ is observable, we can conduct parameter estimation using the slow system
\eqref{sloweqn-app}. In fact, this unknown parameter $a$ is carried over to the slow system \eqref{sloweqn-app} which is now rewritten as
\begin{eqnarray} \label{slowsystem}
\dot x &=& Ax + f(x, \hat h^\eps( x, \theta_t\omega)+\sigma\eta(\psi_\varepsilon\omega), a), \;\; x \in \R^n,
\end{eqnarray}
for $\eps$ sufficiently small.

Since this slow system is lower dimensional than the original system, this method offers an advantage in computational cost, in addition to a benefit of using only observations on slow variables.
It is often   more feasible to observe slow variables   than fast variables \cite{Bishwal}.

Assume that we have
  observation, $x_{ob}$, on the slow component $x$ only, and let us estimate the system parameter $a$.  An estimator $a_E^S(\eps)$  for $a$  may be
obtained with existing  techniques  as reviewed in our earlier work
\cite{JYangDuan}  or \cite{Bishwal, Ibragimov} and references
therein. Here the subscript $S$ indicates that the parameter
estimation  is conducted on the slow system \eqref{slowsystem}. For
example, $a_E^S(\eps)$ may be obtained by minimizing the objective
function $F(a)  \triangleq \EX  ||x_{ob} - x ||^2_{\R^n}$.

We then compare this estimator $a_E^S(\eps)$ with the estimator, $a_E(\eps)$ (without
superscript $S$), based on the original  slow-fast system
\eqref{Slow-Equation-Stoch000}-\eqref{Fast-Equation-Stoch000},
when observations $x_{ob}, y_{ob}$ are available for both components
$x, y$. For example, $a_E(\eps)$ may be obtained by minimizing the objective function
$ \mathbb{F} (a) \triangleq  \EX [ ||x_{ob} - x ||^2_{\R^n} + ||y_{ob} - y ||^2_{\R^n}] $.

A stochastic Nelder-Mead method is used to minimize the objective functions.
 In the next section, we demonstrate this method with an example.

\section{An example}
In this section, we demonstrate our parameter estimation method
based on random slow manifolds by a simple example.
\begin{example}
Consider  a slow-fast stochastic system
\begin{eqnarray}\label{slow-eps-SDE}
 \dot x&=& 0.001x- a xy, \quad x(0)=x_0 \in \R,\\\label{fast-eps-SDE}
 \dot y &=& \frac1{\eps} (-y + \frac{1}{600} x^2) + \frac{\sigma}{\sqrt{\eps}}
 \dot{W}_t,\quad y(0)=y_0 \in \R,
\end{eqnarray}
 where $a $ is a real unknown positive parameter, $\eps$ is a small positive scale separation constant,
 $\sigma$ is a constant   noise intensity, and $W_t$ is a scalar Wiener process.

 See Figure \ref{smallchange001} for a phase portrait of the corresponding deterministic system ($\sigma=0$).

 \begin{figure}[h]
  \centering
  \epsfig{file=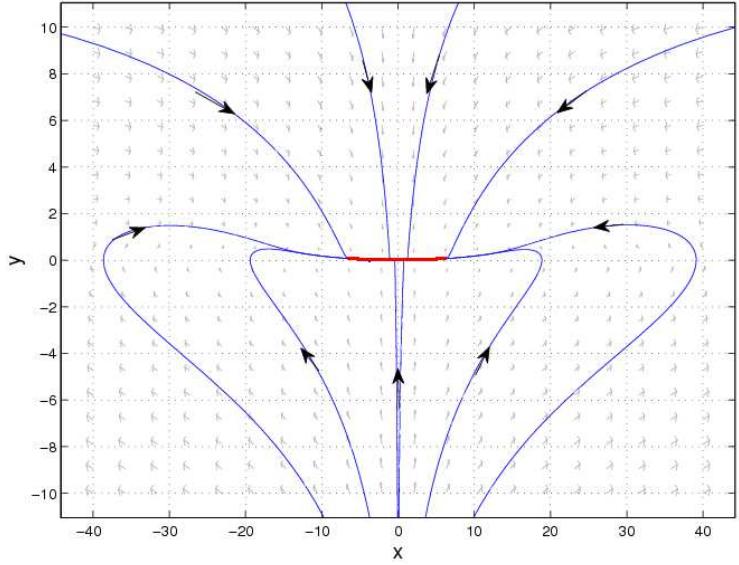,width=0.6\linewidth}
  \caption{Deterministic dynamics -- Phase portrait for   $\dot x=0.001x-xy$, $\dot y=\frac1{\eps} (-y + \frac1{600} x^2)$ with $\eps=\frac16$:
  The   global attractor is clearly seen  (Red or thick curve within $-7<x<7$ and near the $x-$axis).} \label{smallchange001}
 \end{figure}

In this system,  the nonlinear terms are not global Lipschitz. But
if additionally it has  an absorbing set, we can   cut-off the
nonlinearities without affecting the long time, slow dynamics
(almost surely). Indeed, for arbitrary constants $M$ and $K$, we
have
\begin{eqnarray*}
dx^2&=&2x\,dx=(0.002x^2-2ax^2y)\,dt,\\
d(My-K)^2&=&2(My-K)M\,
dy+M^2\,d[y,y]\\{}&=&\big(-\frac{2M^2}{\eps}y^2+\frac{M^2}{300\eps}x^2y+
\frac{2MK}{\eps}y-\frac{MK}{300\eps}x^2+\frac{M^2\sigma^2}{\eps}\big)dt+2M(My-K)\frac{\sigma}{\sqrt{\eps}}\,dW_t.
\end{eqnarray*}
Therefore,
\begin{eqnarray}
\lefteqn{d\big(\frac{M^2}{300\eps}x^2+2a(My-K)^2\big)}\nonumber\\&=&-\frac1{\eps}\big(\frac{M^2}{300\eps}x^2+2a(My-K)^2\big)\,dt-\frac{2aM^2}{\eps}y^2\,dt
+\frac{0.002M^2-2aMK+M^2/\eps}{300\eps}x^2\,dt
\nonumber\\{}&&+\frac{2aK^2+2aM^2\sigma^2}{\eps}\,dt +
4aM(My-K)\frac{\sigma}{\sqrt{\eps}}\,dW_t.
\end{eqnarray}
Taking $M=\eps$, $K=\frac1{a\eps}$, then we see that
\[\frac{0.002M^2-2aMK+M^2/\eps}{300\eps}=\frac{0.002\eps^2-2+\eps}{300\eps}<0,  \text{  for  small }\eps,\]
and \[\frac{2aK^2+2aM^2\sigma^2}{\eps}=\frac2{a\eps^3}+
2a\eps\sigma^2.\] Thus
\begin{eqnarray}
\frac{d}{dt}\mathbb{E}\big(\frac{M^2}{300\eps}x^2+2aM^2(y-K/M)^2\big)\leq
-\frac1{\eps}\mathbb{E}\big(\frac{M^2}{300\eps}x^2+2aM^2(y-K/M)^2\big)+
\frac2{a\eps^3}+ 2a\eps\sigma^2.
\end{eqnarray}
By the Gronwall inequality, we conclude that
\begin{eqnarray}
\lefteqn{\mathbb{E}\big(\frac{\eps}{300}x^2+2a\eps^2(y-\frac1{a\eps^2})^2\big)}\nonumber\\{}&\leq&
\mathbb{E}\big(\frac{\eps}{300}x_0^2+2a\eps^2(y_0-\frac1{a\eps^2})^2\big)e^{-\frac{t}{\eps}}+
(\frac2{a\eps^3}+
2a\eps\sigma^2)\varepsilon(1-e^{-\frac{t}{\eps}})\nonumber\\{}&\leq&
\big(\frac{\eps}{300}x_0^2+2a\eps^2(y_0-\frac1{a\eps^2})^2\big)e^{-\frac{t}{\eps}}+
(\frac2{a\eps^2}+ 2a\eps^2\sigma^2).
\end{eqnarray}
This means, for fixed $\eps$ and $a>0$, the dynamics of the system
\eqref{slow-eps-SDE}-\eqref{fast-eps-SDE} will eventually stay in an
ellipse (almost surely), i.e. there is a random absorbing set. We
can then cut-off the nonlinearities outside this absorbing set to
obtain a modified system which has,   almost surely, the same long
time, slow dynamics as the original system \cite{Kan2013}. In the
following calculations, we actually have omitted this cut-off
procedure for simplicity.
\bigskip

By the random transformation \eqref{random transformation}, SDEs
system \eqref{slow-eps-SDE}-\eqref{fast-eps-SDE} are converted
into the following  system
\begin{eqnarray}\label{slow-eps-RDE}
 \dot X&=& 0.001X- a X(Y+\sigma\eta^\varepsilon(\theta_t\omega)) ,\\\label{fast-eps-RDE}
 \dot Y &=& -\frac1{\eps} Y + \frac1{\eps} \frac{1}{600} X^2 .
\end{eqnarray}
Therefore, there exists an $\tilde h^\eps$  satisfying
\begin{eqnarray} \label{slow-example}
 \tilde h^\eps(\xi, \omega)= \frac1{\eps}\int_{-\infty}^0
e^{\frac{s}{\eps}}\,\frac{X^2}{600}\,ds ,
\end{eqnarray}
whose graph is a random slow manifold  for the random system
\eqref{slow-eps-RDE}-\eqref{fast-eps-RDE}. In fact, $\tilde h^\eps(\xi,
\omega)$ has an approximation $\hat h^\eps(\xi, \omega)$ (with error
$\mathcal{O}(\varepsilon^2)$),
\begin{eqnarray}\label{slow-manifold}
\hat h^\eps(\xi, \omega)= \frac{\xi^2}{600} + \varepsilon
\big(-\frac{\xi^2}{300} 0.001 + \frac{\xi^4}{180000}a - \frac{\xi^2
}{300} a \sigma \int_{-\infty}^0 s e^s\, dW_s\big).
\end{eqnarray}
So the approximated slow   system is
\begin{eqnarray}
 \dot x &=& 0.001x - a x \big(\sigma \eta(\psi_\varepsilon\omega)   +\hat h^\eps(x, \theta_t\omega)\big). \label{sloweqn2}
\end{eqnarray}with $\eta(\psi_\varepsilon\omega)=\int_{-\infty}^0 e^s\,
dW_s(\psi_\varepsilon\omega)$. We will illustrate that the parameter estimator for $a$ based on this low dimensional, reduced system \eqref{sloweqn2} is a good approximation of the parameter estimator based on the original system \eqref{slow-eps-SDE}-\eqref{fast-eps-SDE}.

\paragraph{Slow-fast system and random slow manifold}
The  random slow manifold is the graph of $h^\varepsilon(\xi,
\omega)=\sigma\eta(\psi_\varepsilon\omega)+ \tilde h^\eps(\xi,
\omega)$, where $\tilde h^\varepsilon$ is as in
\eqref{slow-example}. It is a curve (depending on random samples).
The random orbits of the system
\eqref{slow-eps-SDE}-\eqref{fast-eps-SDE}  approaches to this curve
exponentially fast. Figures \ref{figure
SDEs-sms001e001}--\ref{figure SDEs-sms005e001} show some orbits of
the slow-fast system \eqref{slow-eps-SDE}-\eqref{fast-eps-SDE} and
the approximate random slow manifold $\mathfrak{h}^\varepsilon(\xi,
\omega)=\sigma\eta(\psi_\varepsilon\omega)+ \hat h^\eps(\xi,
\omega)$, where $\hat h^\varepsilon$ is as in \eqref{slow-manifold},
with different $\sigma$ and $\eps$ values.

The orbits of the   system \eqref{slow-eps-SDE}-\eqref{fast-eps-SDE}
decay quickly to the random slow manifold.

Figure \ref{figure SDEs-slma-eps001} shows several samples or realizations of
the random slow manifold   with different $a$ values.

\begin{figure}[H]
\includegraphics[height=6.8cm]{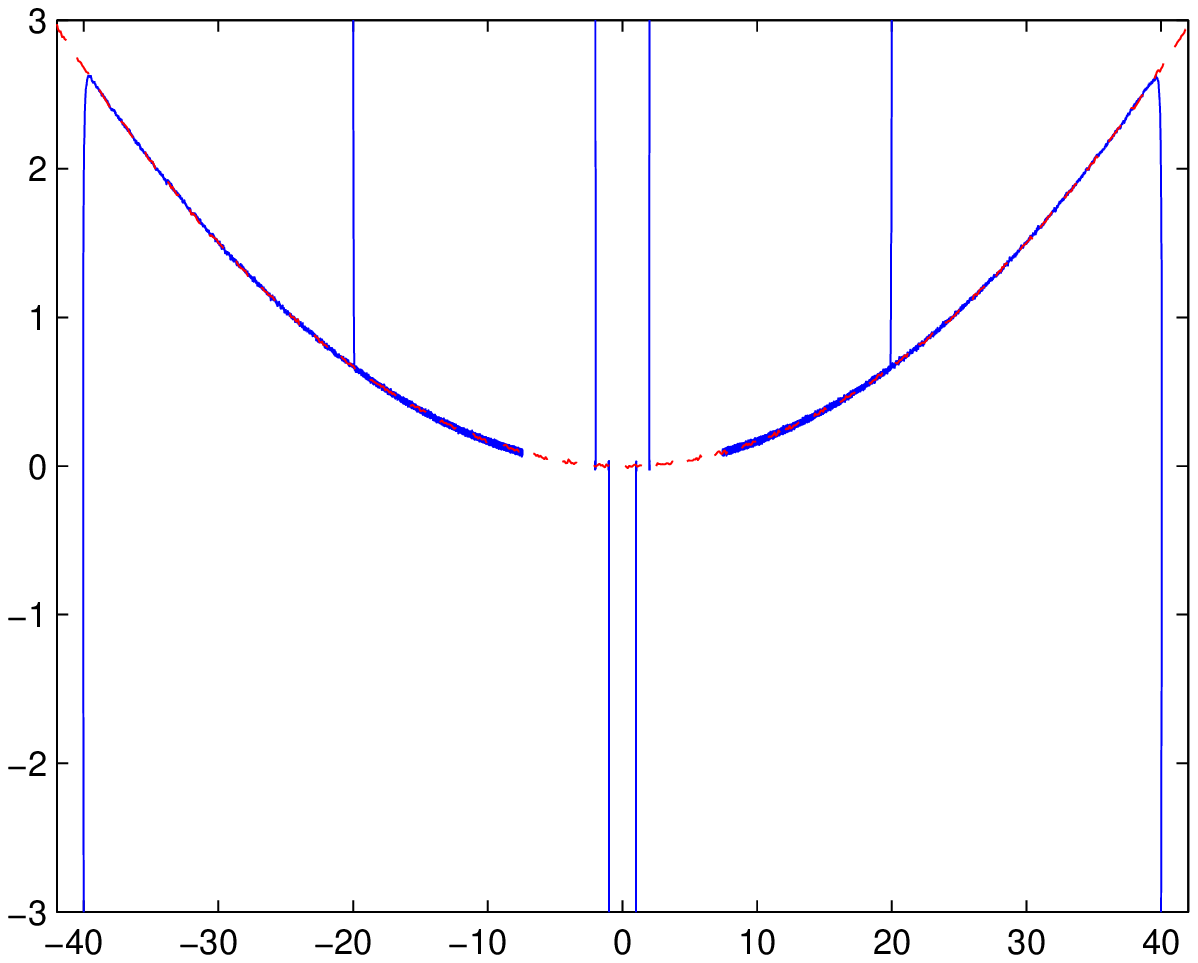}
\includegraphics[height=6.8cm]{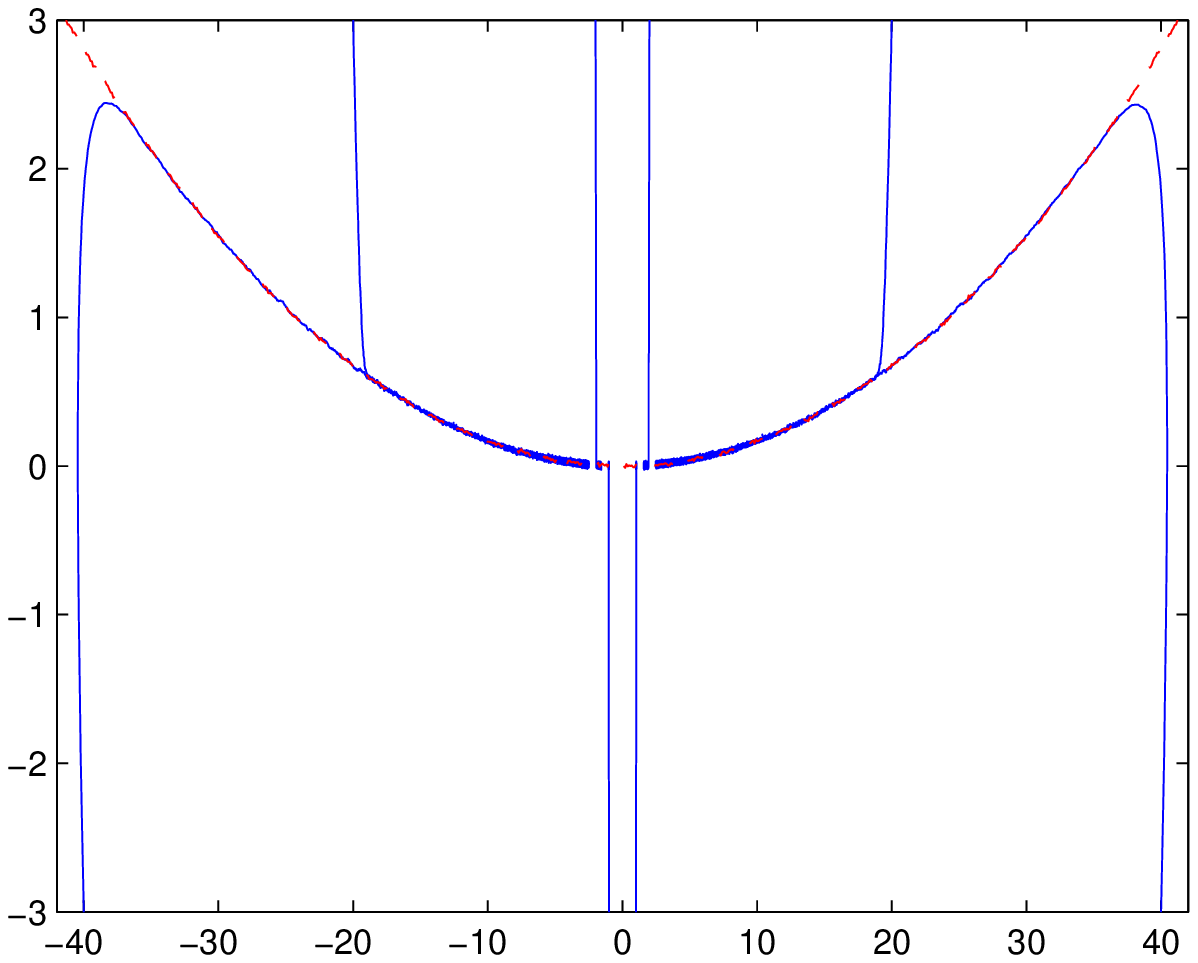}
\caption{Orbits of slow-fast system
\eqref{slow-eps-SDE}-\eqref{fast-eps-SDE} (blue curves) and its slow
manifold expansion $\mathfrak{h}^\varepsilon(\xi,
\omega)=\sigma\eta(\psi_\varepsilon\omega)+ \hat h^\eps(\xi,
\omega)$ (red curve), where $\hat h^\varepsilon$ is as in
\eqref{slow-manifold}, with $\sigma=0.01$ and $\eps=0.01$: $a=0.1$
(left ) and $a=1$ (right).} \label{figure SDEs-sms001e001}
\end{figure}

\begin{figure}[H]
\includegraphics[height=6.8cm]{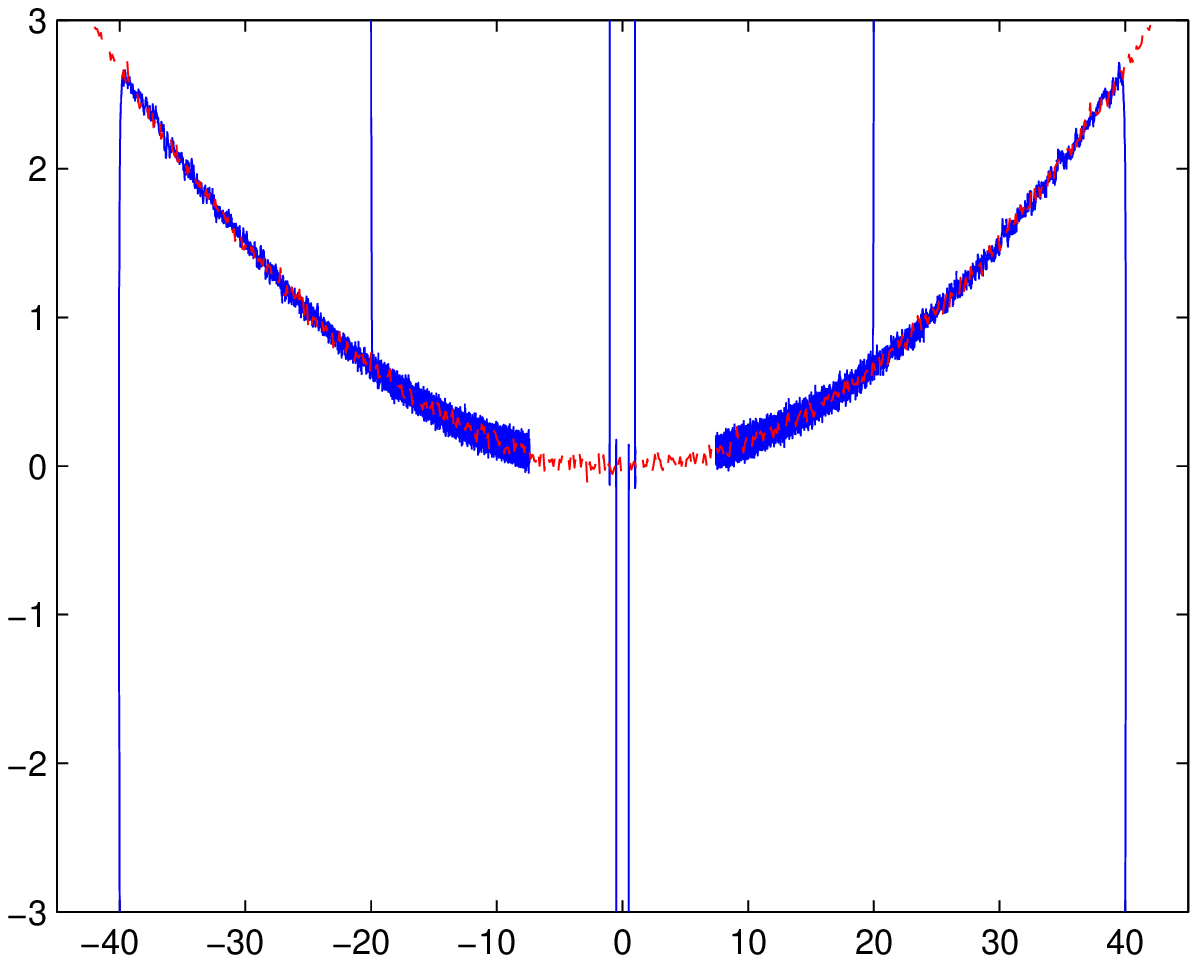}
\includegraphics[height=6.8cm]{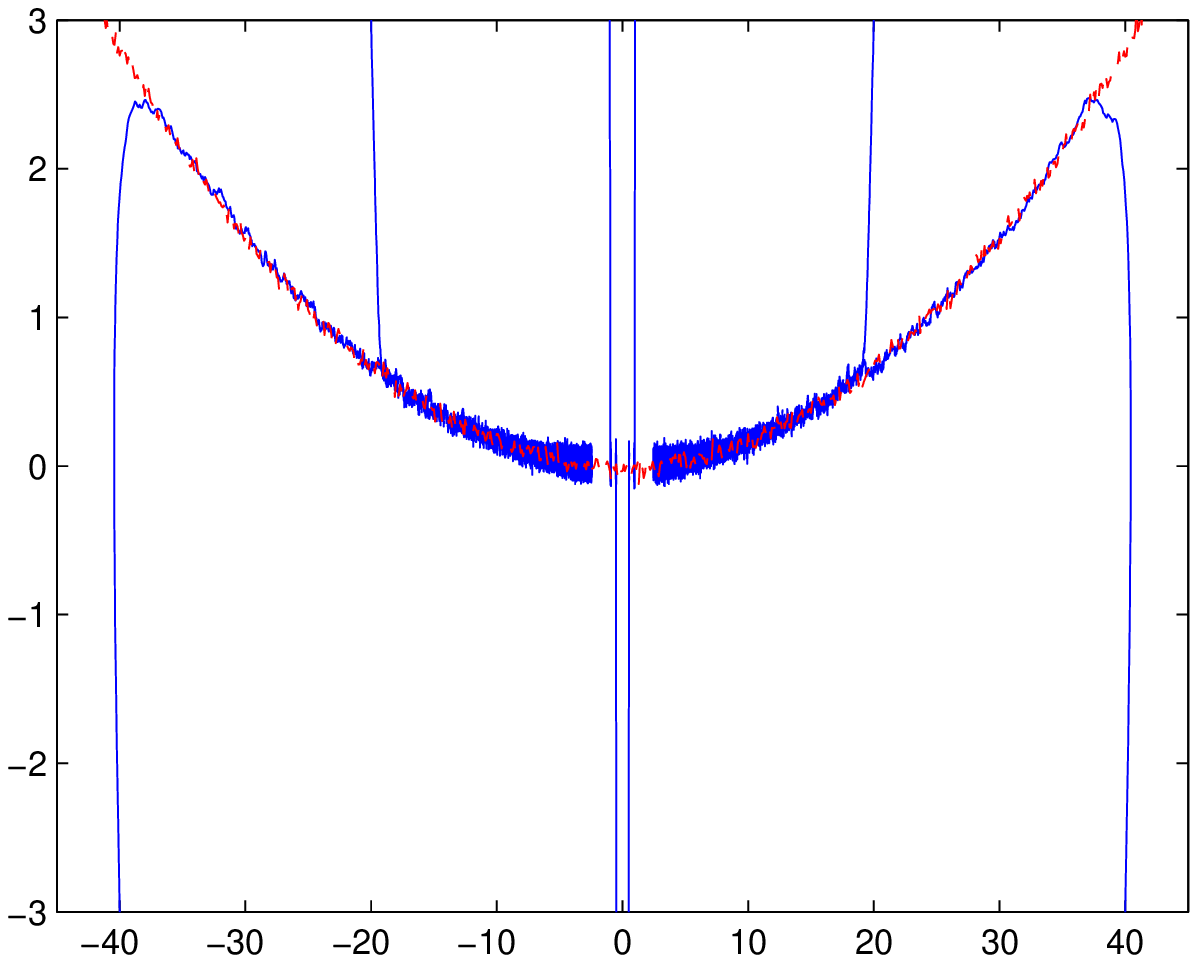}
\caption{Orbits of slow-fast system
\eqref{slow-eps-SDE}-\eqref{fast-eps-SDE} (blue curves) and its slow
manifold expansion $\mathfrak{h}^\varepsilon(\xi,
\omega)=\sigma\eta(\psi_\varepsilon\omega)+ \hat h^\eps(\xi,
\omega)$ (red curve), where $\hat h^\varepsilon$ is as in
\eqref{slow-manifold}, with $\sigma=0.05$ and $\eps=0.01$: $a=0.1$
(left ) and $a=1$ (right).}\label{figure SDEs-sms005e001}
\end{figure}

\begin{figure}[H]
\includegraphics[height=6.8cm]{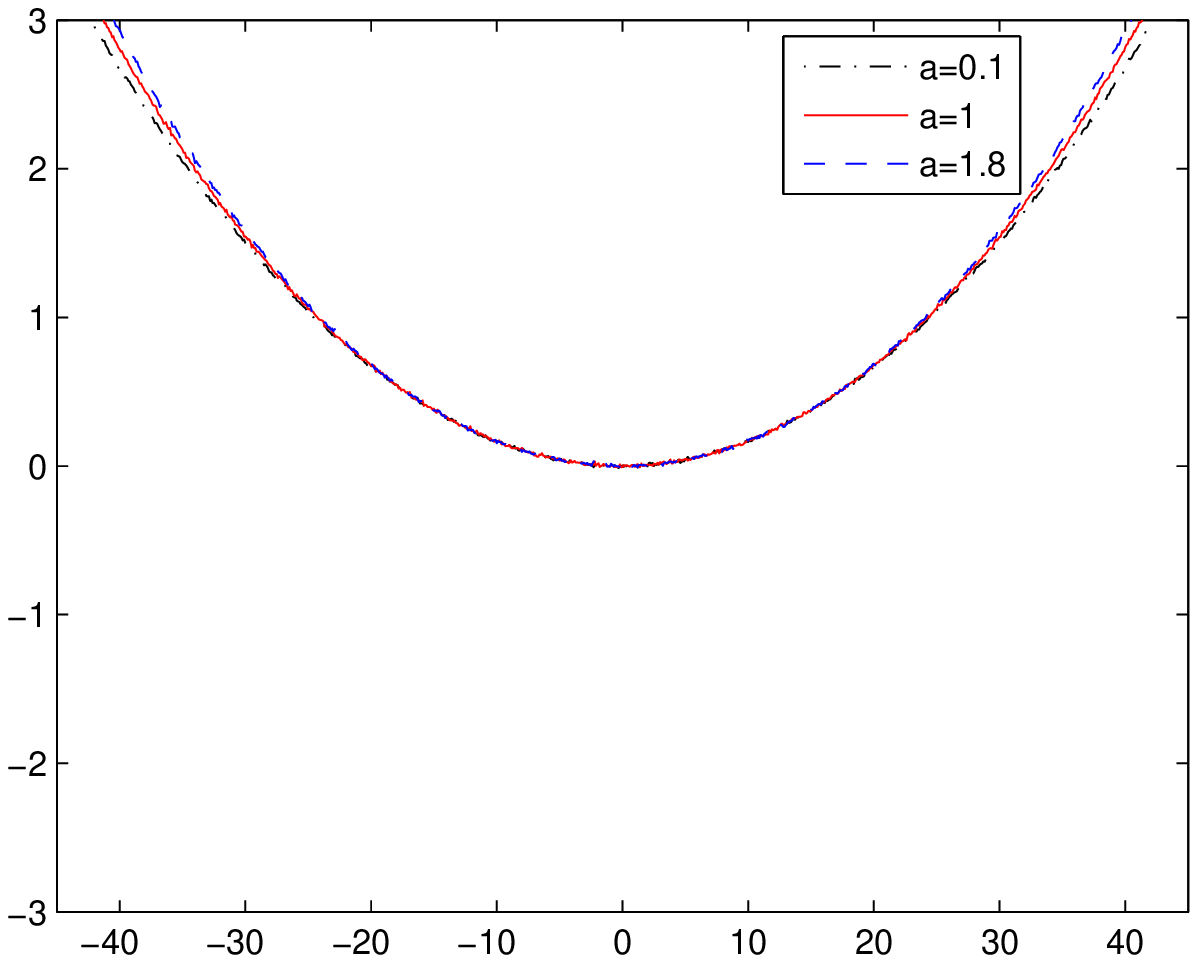}
\includegraphics[height=6.8cm]{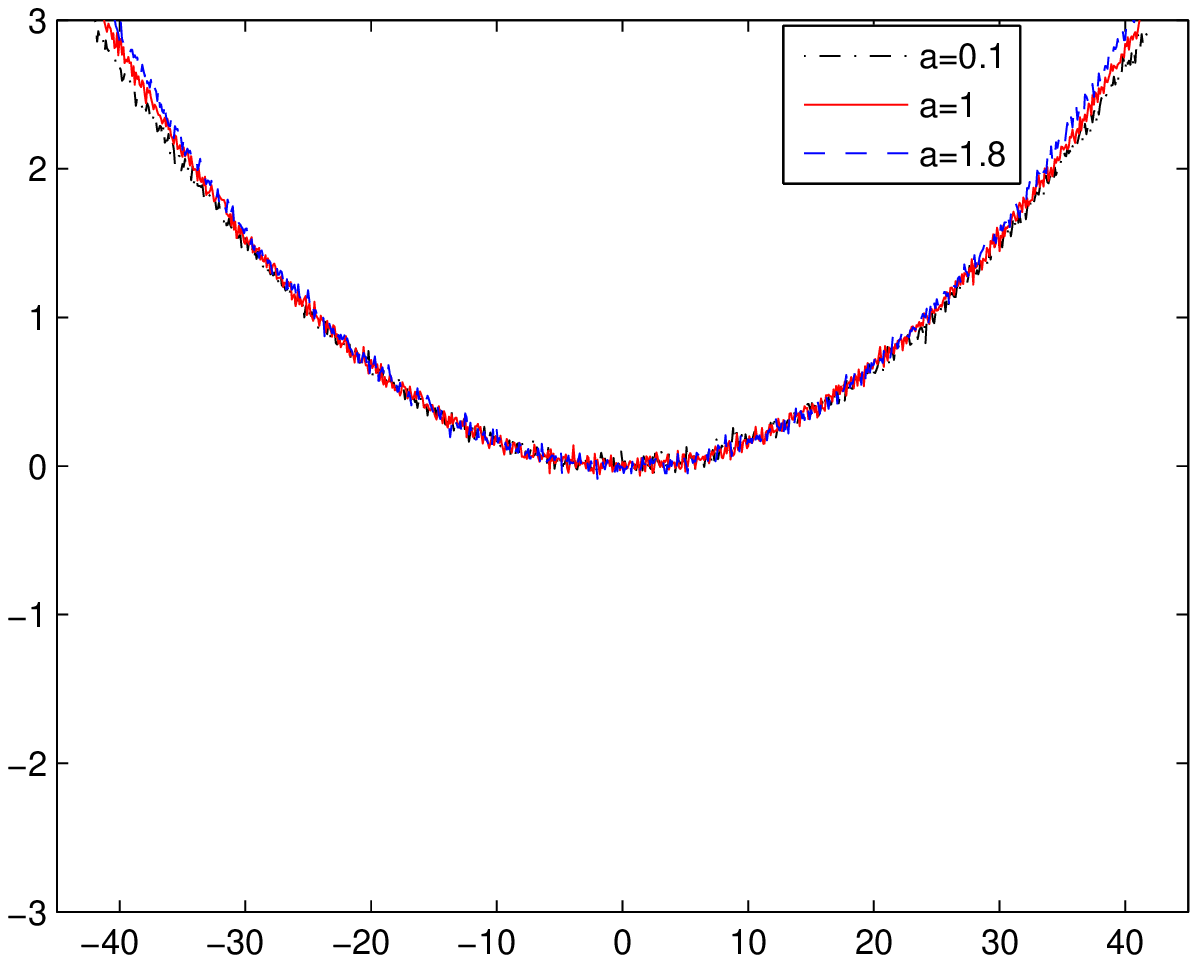}
\caption{Random slow manifold expansion
$\mathfrak{h}^\varepsilon(\xi,
\omega)=\sigma\eta(\psi_\varepsilon\omega)+ \hat h^\eps(\xi,
\omega)$, where $\hat h^\varepsilon$ is as in \eqref{slow-manifold}
with $\eps=0.01$: $\sigma=0.01$ ( left ) and $\sigma=0.05$
(right).}\label{figure SDEs-slma-eps001}
\end{figure}

\paragraph{Nelder-Mead method and stochastic Nelder-Mead method}
 The deterministic Nelder-Mead method (NM) is a geometric search method to find a minimizer of
  an objective function  $F(a)$. Starting from an initial guess point, it   generates a new point
(reflection point, expand point, inside/outside contraction point or
shrink point) by comparing function values, and thus get a
better and better estimator until the smallest objective function
value in this iteration reaches the termination value (prescribed error tolerance). The
algorithm of the NM method and its improvements have been widely
 studied and utilized   \cite{Barati, Nocedal, Pham}. The
method has its advantage that the objective function need not to be
differentiable, and it can thus be used in various applications.  As
noted in \cite{Chang, Price}, the Nelder-Mead method is a
  widely used heuristic algorithm. Only very limited convergence results
exist for a   class of low dimensional (one or two dimensions)
problems, such as in the case when the objective function $F$ is strictly
convex \cite{Lagarias}. For the numerical simulation, with   analytical objective
functions, Matlab function \emph{fminsearch} can be
used to find a minimizer.

However, when dealing with problems with noise, NM method has the
disadvantage \cite{Barton, Chang} that it lacks an effective sample
size scheme for controlling noise, as shrinking steps are sensitive
to the noise in the objective function values and then may lead to
the search in a wrong direction. In fact   an analytical and
empirical evidence is known for the false convergence \cite{Barton} on stochastic function.
So we use the stochastic Nelder-Mead simplex method (SNM)
  \cite{Chang} to    mitigate the possible mistakes   in  the
stochastic setting. The newly developed Adaptive
Random Search   in \cite{Chang}  consists of a local search and a global search. It
  generates a new point   and new
objective function $\hat F=\sum\limits_{i=1}^{N(k)}F_i/N(k)$ in the
$k-$th iteration with increasing number of the sample size scheme
$N(k)$.  A proper choice for $N(k)$  is $[\sqrt{k}]$,
with $[c]$ the largest integer not bigger than  $c$. Here
$\sum\limits_{i=1}^{N(k)}F_i$ is the sum of $N(k)$ objective function values for $F$. SNM   leads to the convergence of $F$ at search points
to $\min\limits_{a} F(a)$ (and thus we obtain a  minimizer $a^*$), with probability one.

\paragraph{Parameter estimation}
To illustrate our method for parameter estimation on the random slow manifold, we fix    $\eps=0.01$ and
$\sigma=0.01$ in the following numerical experiments.

We want to estimate the parameter $a$ by using both the original
slow-fast system and the   slow     system, in order to demonstrate that the slow system is appropriate for parameter estimation, when $\eps$ is sufficiently small.

Step 1: Generate observations

Take the true value $a=0.1$ (say) and numerically solve
\eqref{slow-eps-SDE}-\eqref{fast-eps-SDE} with an initial condition
$(x_0, y_0)$ to get $J$ samples of observational data $(x^{ij}_{ob},
y^{ij}_{ob})$, $j=1, \cdots, J$ at time instants $t_i$, $i=1,\cdots,
I$ (save these data).

Step 2:  Estimator  $a_E$ for the original slow-fast system

Take two initial guesses for the unknown system parameter  $a=a_0$ and $a=a_1$ randomly and
  solve the original slow-fast system \eqref{slow-eps-SDE}-\eqref{fast-eps-SDE}, with the same initial
condition $(x_0, y_0)$ and time points $t_i$, $i=1,\cdots, I$. Thus
we obtain    $x^i$ and $y^i$ values which depend on $a$.

Using the stochastic Nelder-Mead Algorithm, we find a parameter
estimator $a_E$ such that the objective function
$F(a)\triangleq\mathbb{E}\sum\limits_{i=1}^{I}\sum\limits_{j=1}^{J}\big((x^i-x^{ij}_{ob})^2
+ (y^{i}-y^{ij}_{ob})^2\big)$ is minimized.

Step 3: Estimator  $a_E^S$ for the   slow  system

Take two initial guesses  $a_0$ and $a_1$ randomly  and solve the
slow system \eqref{sloweqn2} with the same initial condition  $x_0$
at the same time instants $t_i$, $i=1,\cdots, I$. We thus  obtain
$X^i$  which depends on $a_0$ or $a_1$.

By the stochastic Nelder-Mead method as in Step 2, we find  the
parameter estimator $a_E^S$ by minimizing the objective function
$\mathbb{F}(a)
\triangleq~\mathbb{E}\sum\limits_{i=1}^{I}\sum\limits_{j=1}^{J}(X^{i}-~x^{ij}_{ob})^2
$.

\begin{figure}[H]
\includegraphics[height=7cm]{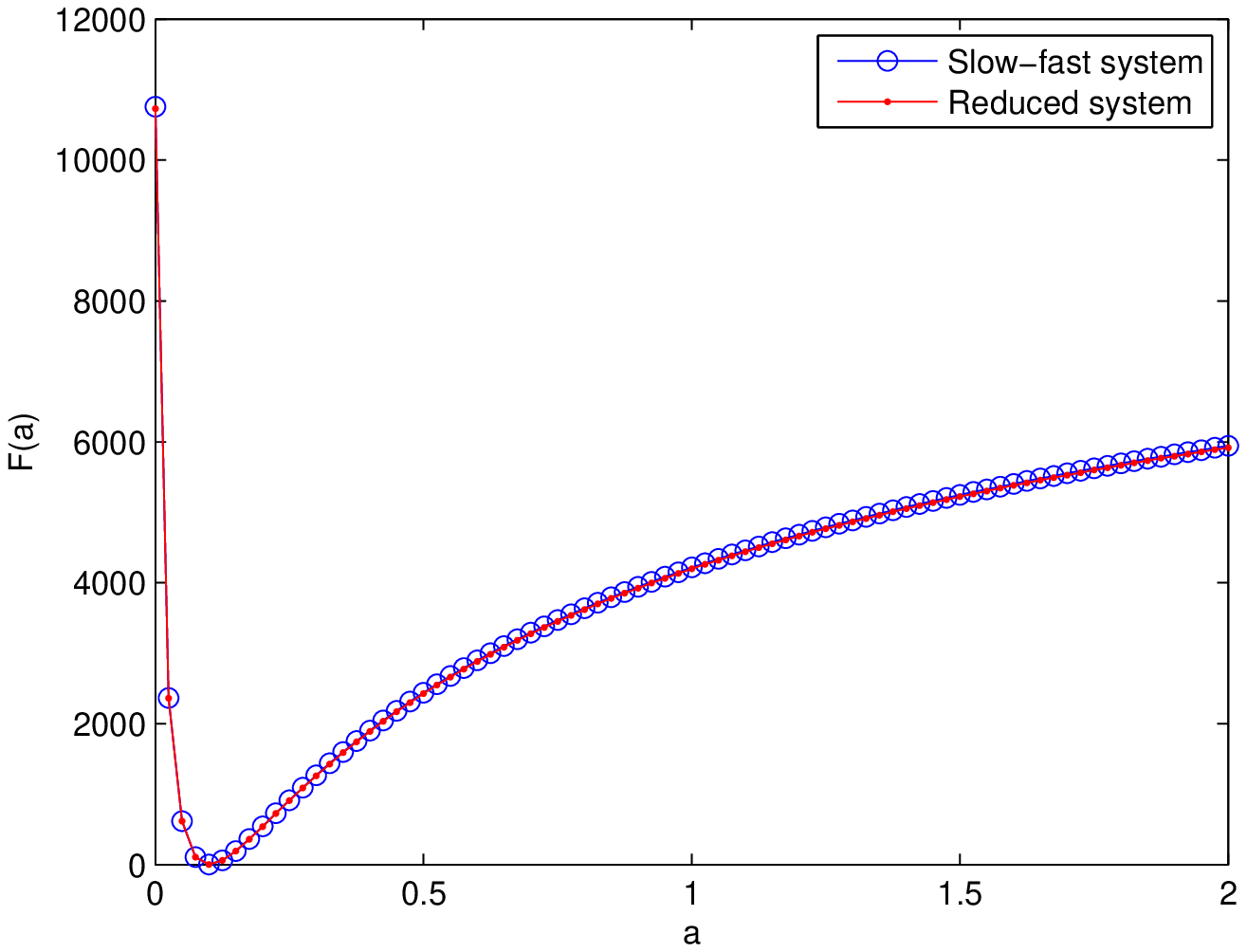}
\includegraphics[height=7cm]{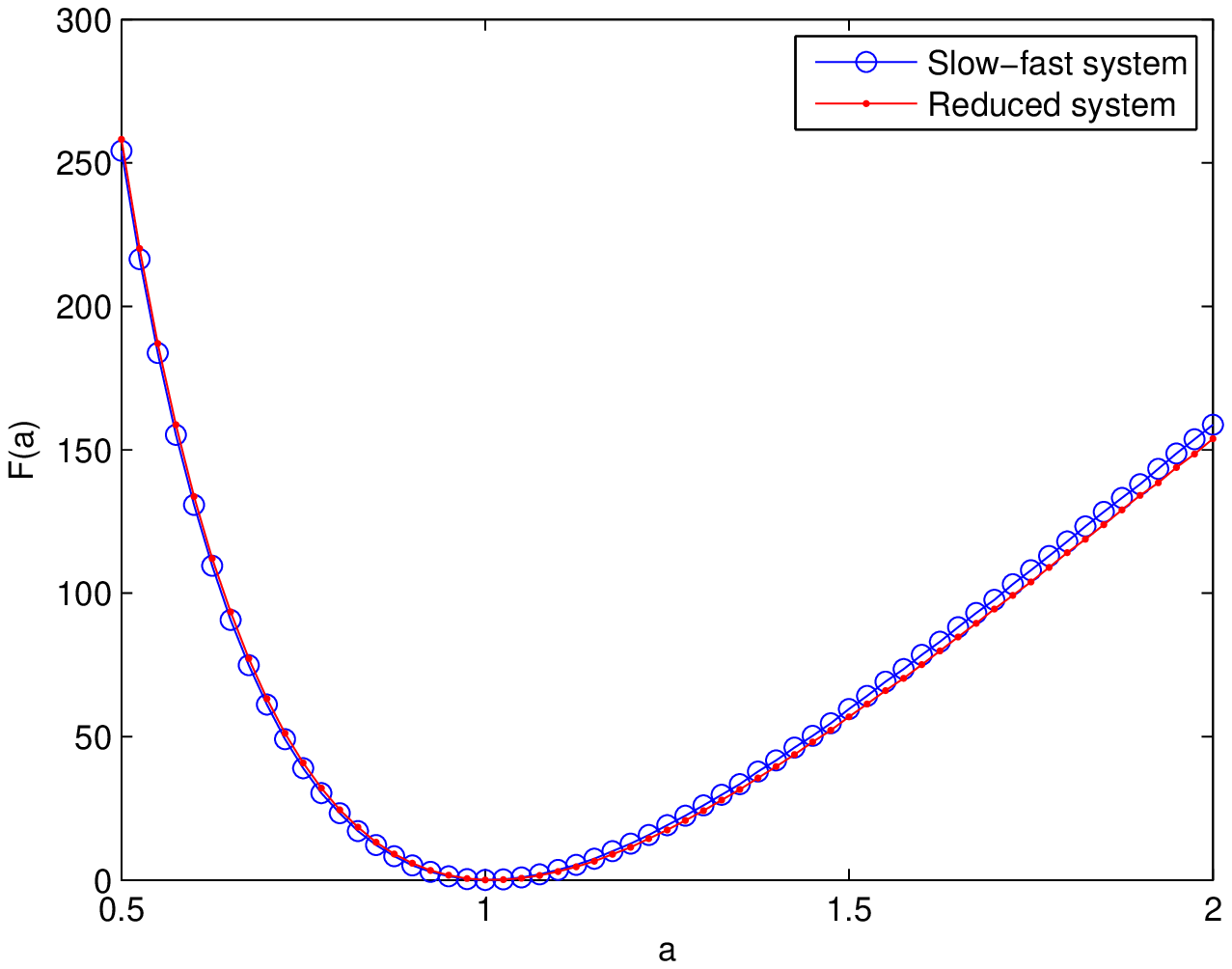}
\caption{Objective function $F(a)$ for slow-fast system
\eqref{slow-eps-SDE}-\eqref{fast-eps-SDE}
 and slow manifold reduced system
\eqref{sloweqn2}, with  $\sigma=0.01$ and $\eps=0.01$: True values
$a=0.1$ (left) and $a=1$ (right).}\label{figure
SDEs-sm-Fs001e001}
\end{figure}

\begin{figure}[H]
\includegraphics[height=7cm]{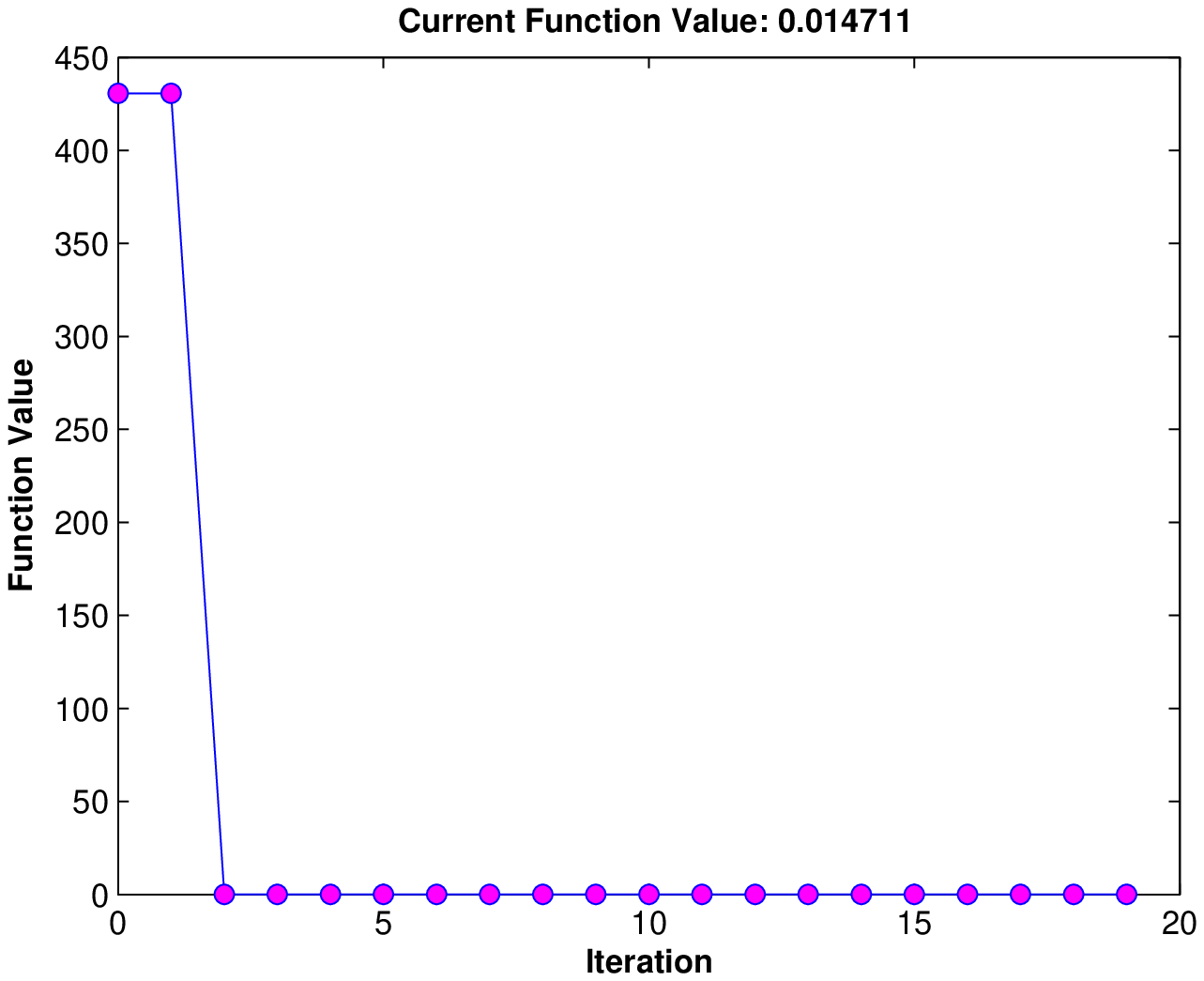}
\includegraphics[height=7cm]{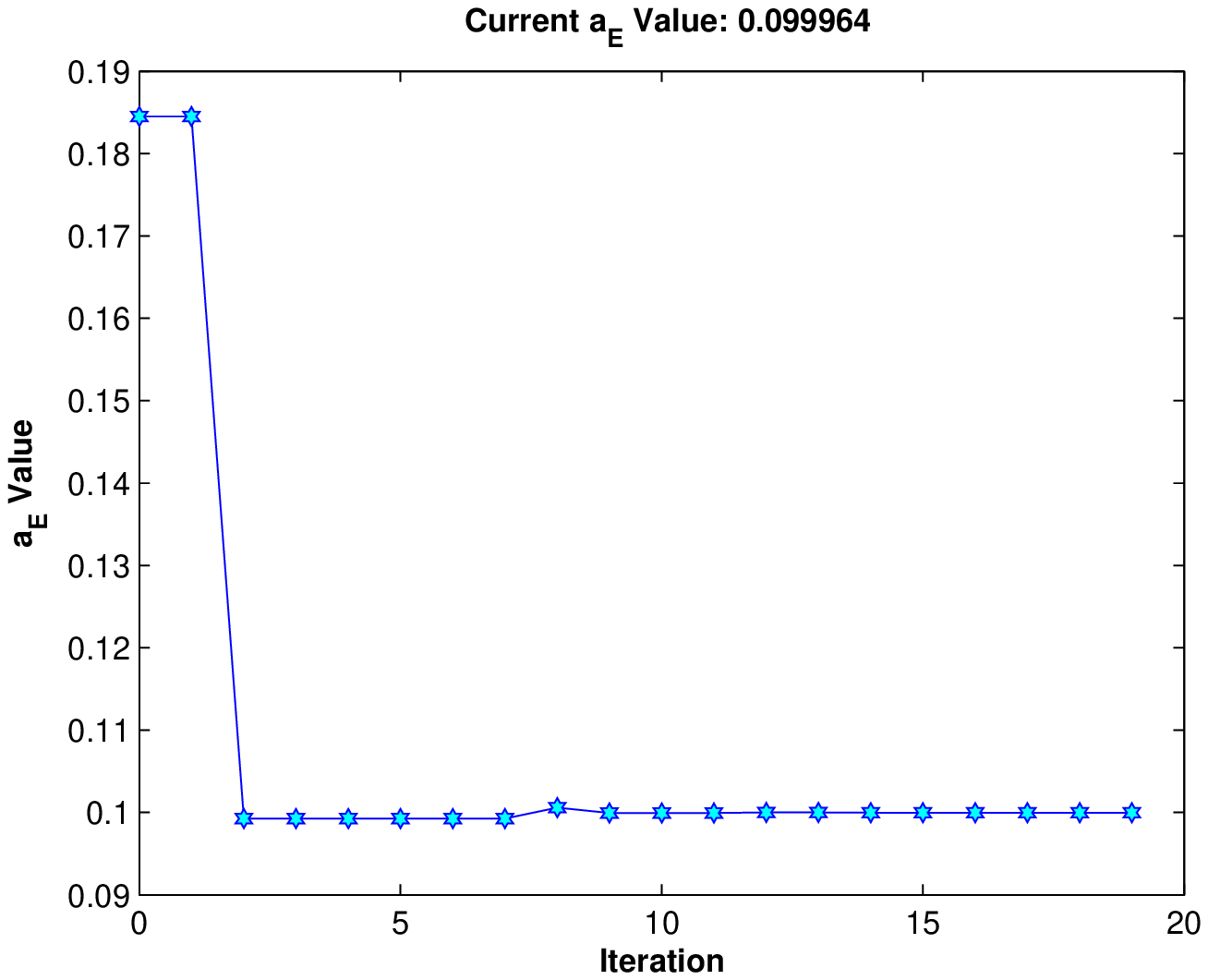}\\
\includegraphics[height=7cm]{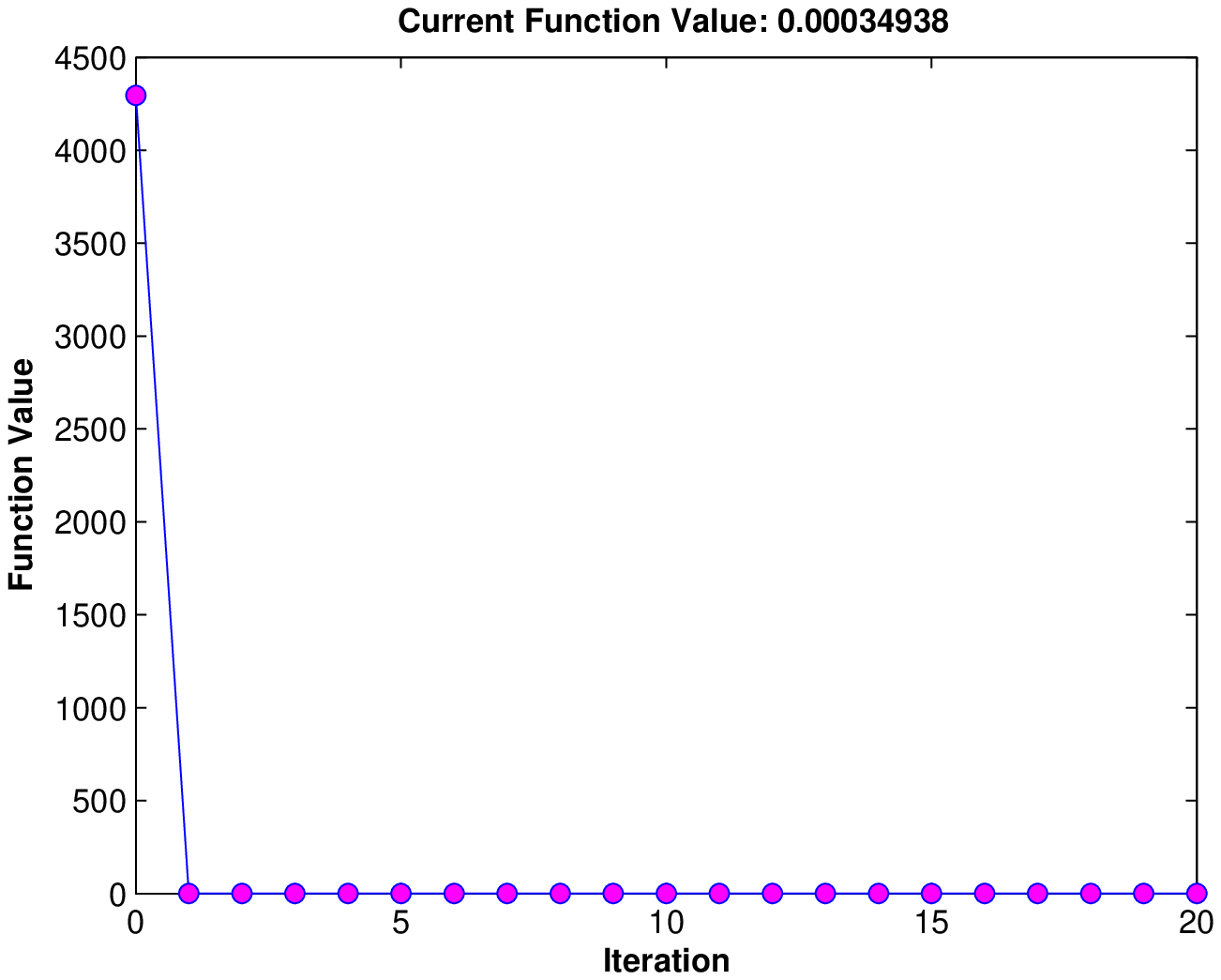}
\includegraphics[height=7cm]{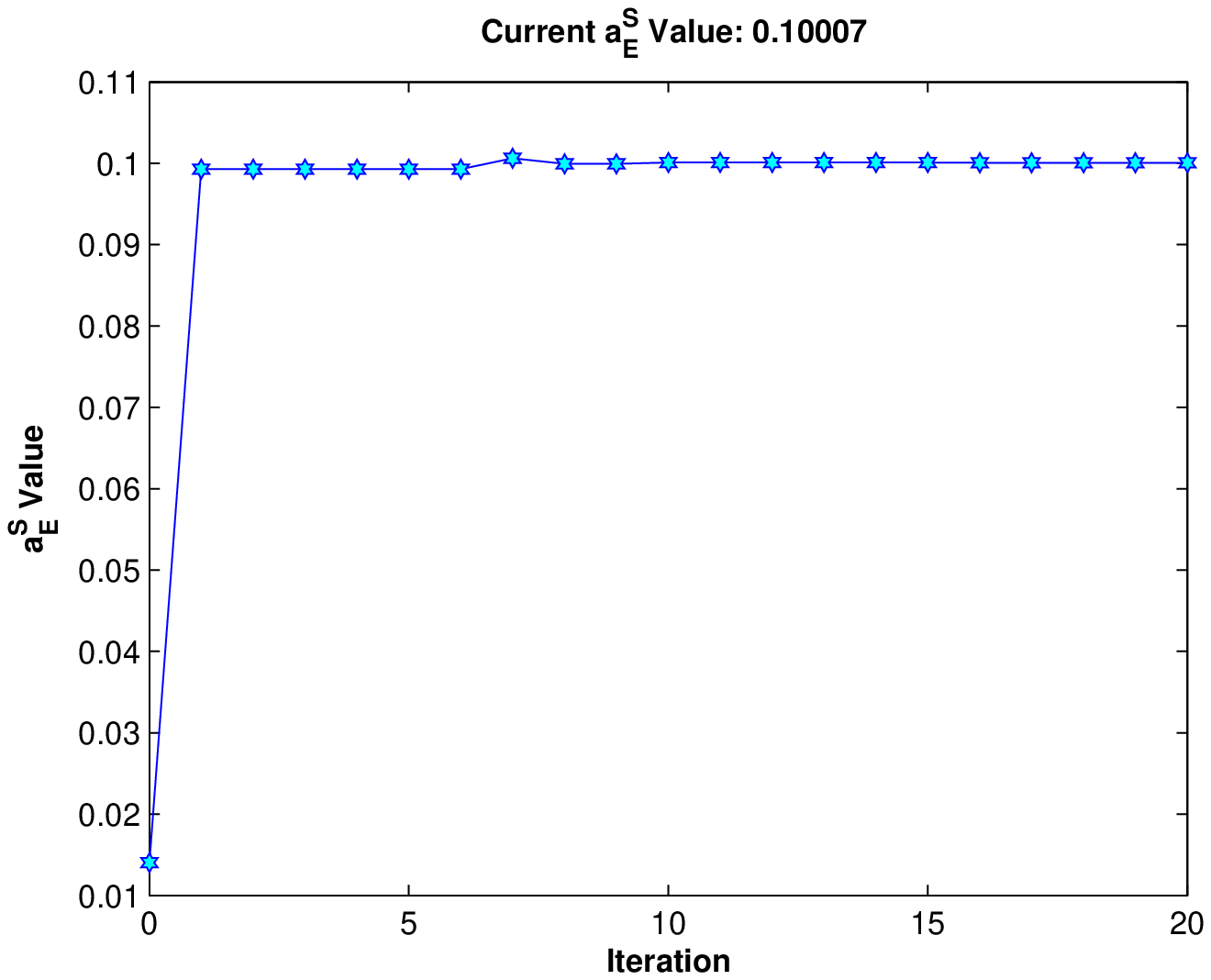}
\caption{Objective function $F(a)$ and estimator $a_E$ (top) for   slow-fast system \eqref{slow-eps-SDE}-\eqref{fast-eps-SDE};
 and objective function $\mathbb{F}(a)$ and estimator $a_E^S$ (bottom) for   slow manifold reduced system
\eqref{sloweqn2}:   $\sigma=0.01$, $\eps=0.01$ and true value
$a=0.1$.}\label{s001e001a01}
\end{figure}

\begin{figure}[H]
\includegraphics[height=7cm]{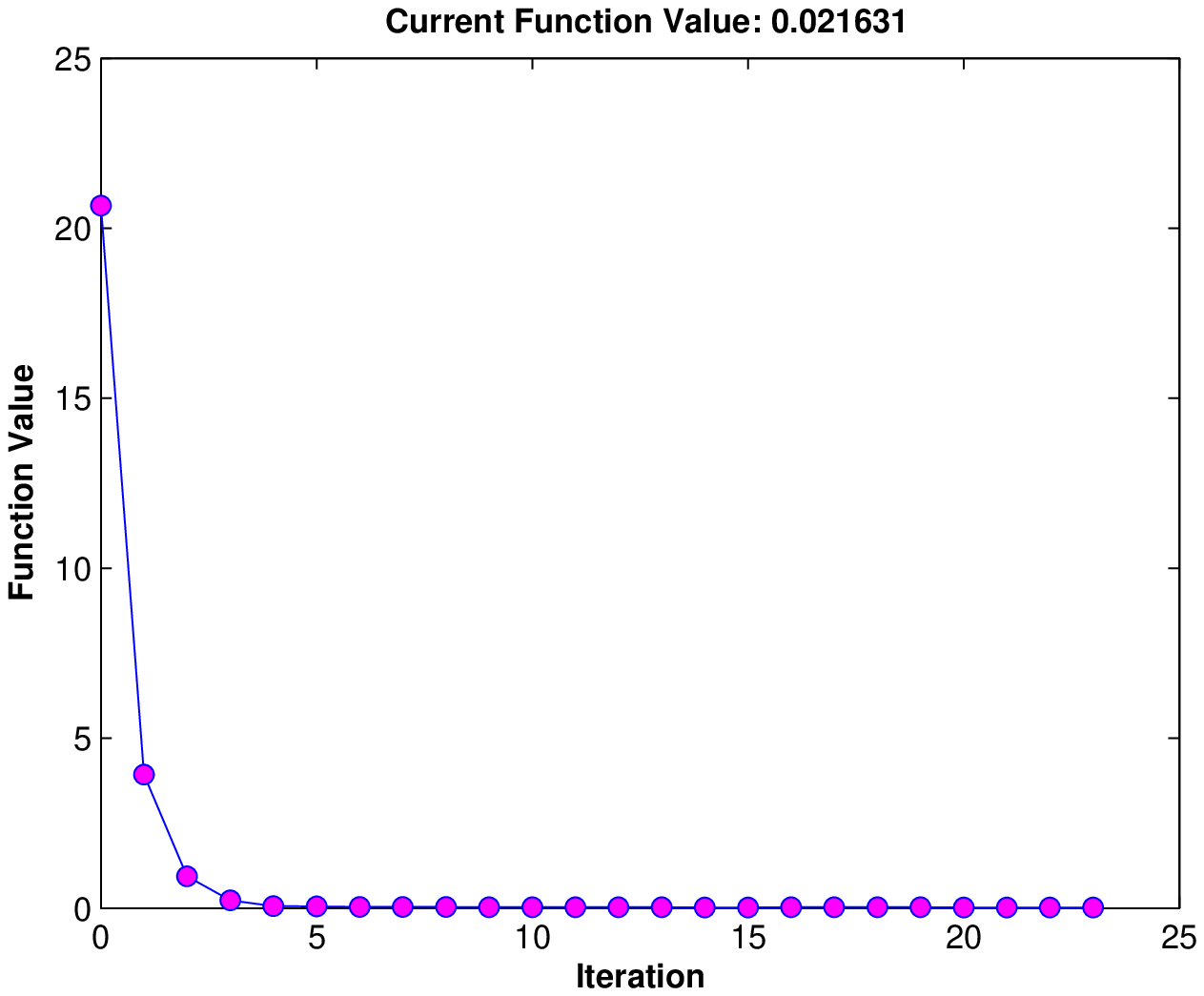}
\includegraphics[height=7cm]{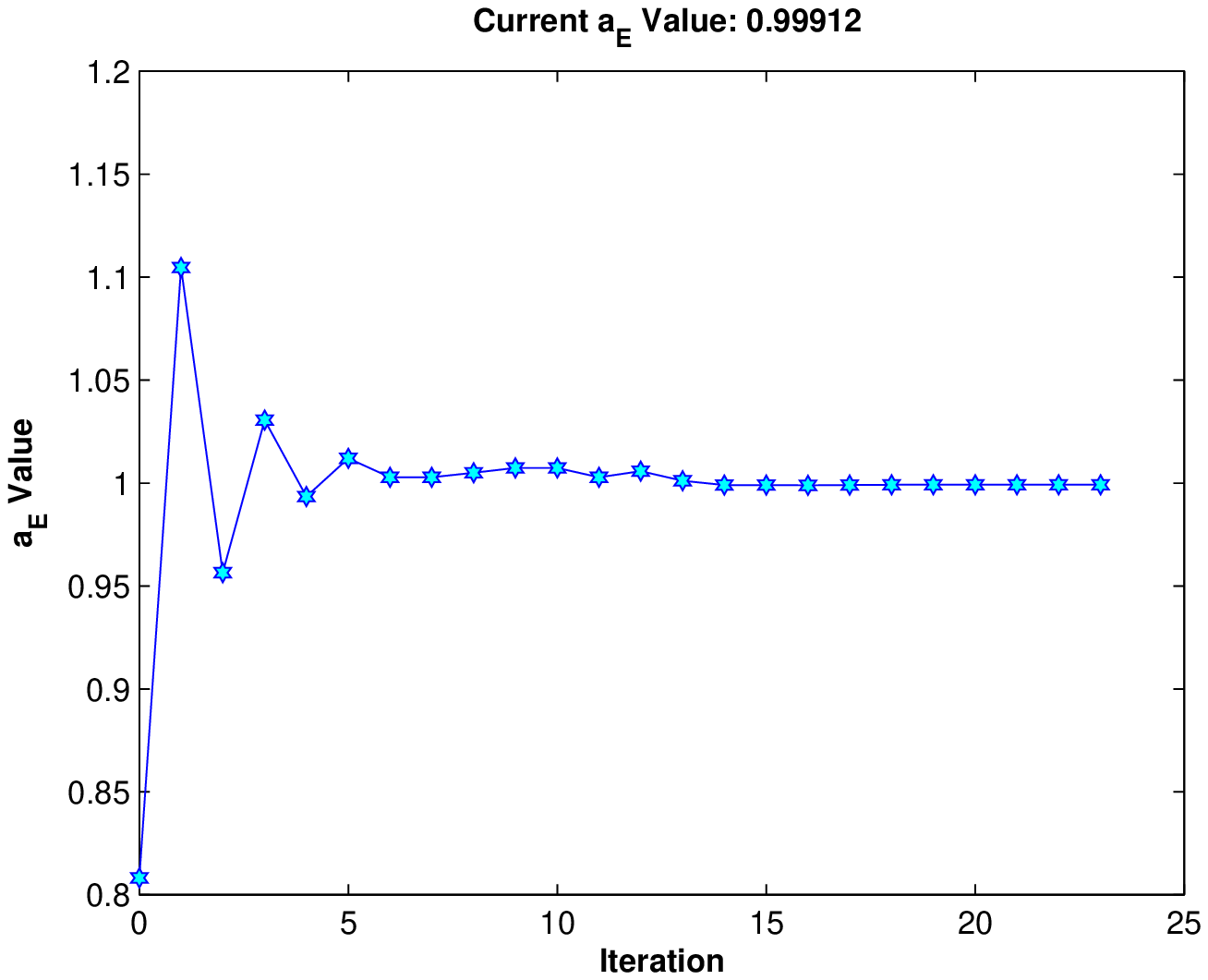}\\
\includegraphics[height=7cm]{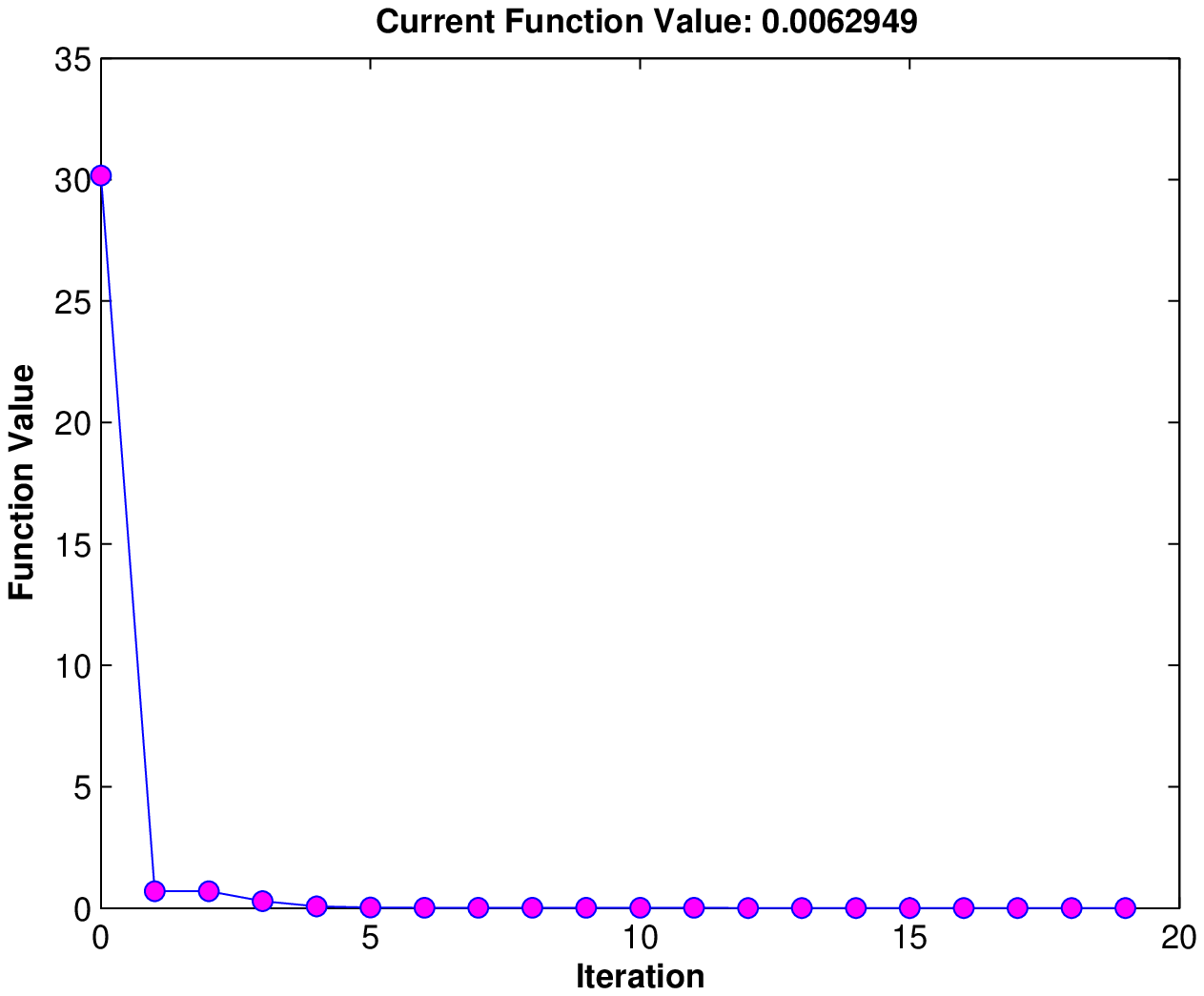}
\includegraphics[height=7cm]{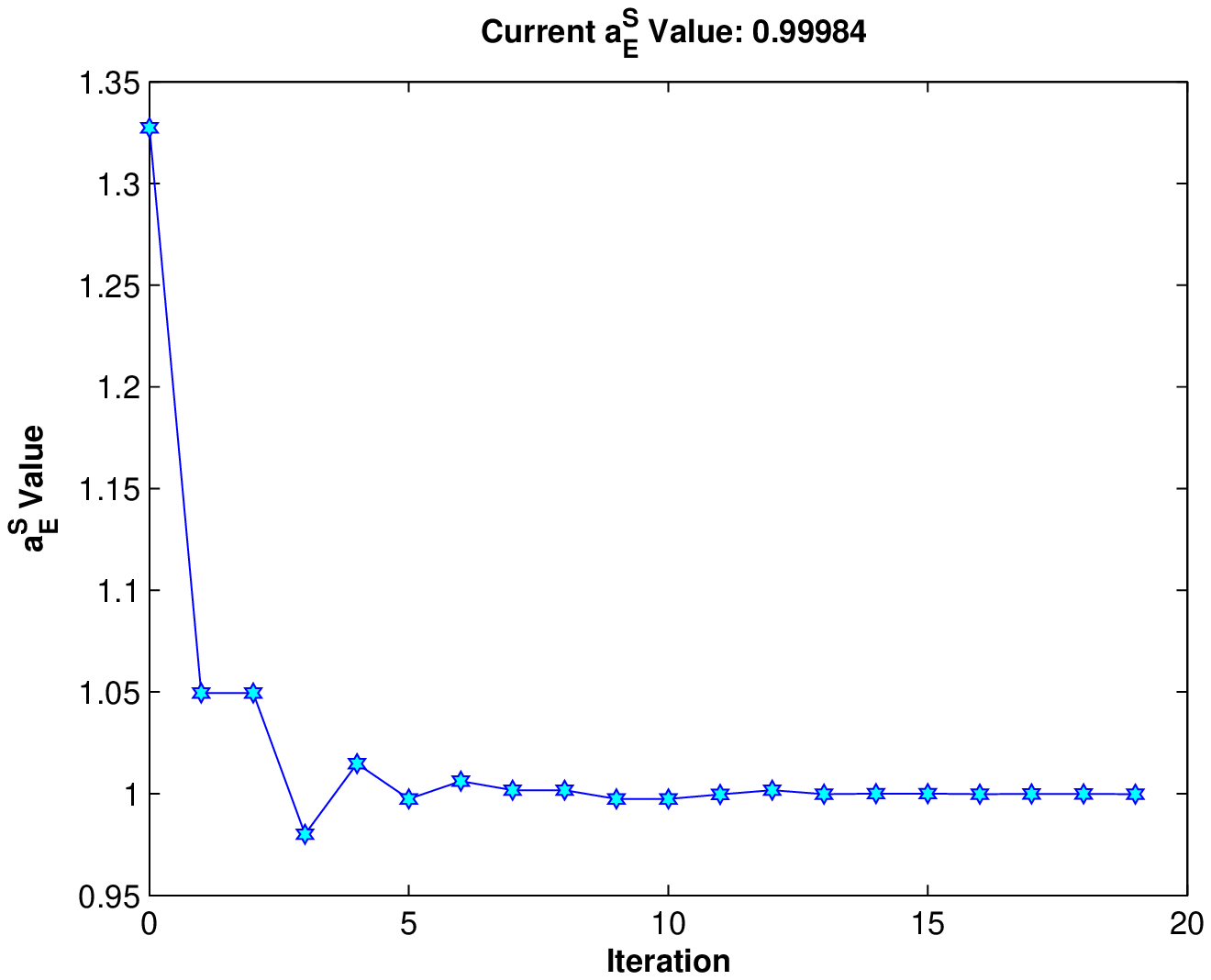}
\caption{Objective function $F(a)$ and estimator $a_E$ (top) for slow-fast system \eqref{slow-eps-SDE}-\eqref{fast-eps-SDE},
 and objective function $\mathbb{F}(a)$ and estimator $a_E^S$ (bottom) for slow manifold reduced system
\eqref{sloweqn2}: $\sigma=0.01$, $\eps=0.01$ and true value
$a=1$.}\label{s001e001a1}
\end{figure}

\paragraph{Numerical experiments}
Figure \ref{figure SDEs-sm-Fs001e001} shows the objective functions
$F(a)$ and $\mathbb{F}(a)$ with true   $a=0.1$ (left) and $a=1$ (right). Here we used $30$ paths to get the expectation in the
objective function.

For $\eps=0.01$ and $\sigma=0.01$,  Figures \ref{s001e001a01} and
\ref{s001e001a1} are objective function values and estimators of
each iteration for slow-fast system
\eqref{slow-eps-SDE}-\eqref{fast-eps-SDE} (top) and
reduced system \eqref{sloweqn2} (bottom)  with $a=0.1$ and
$a=1$, respectively.

We observe that, as proved in \cite{Chang}, the objective function
value $F(a)$   tends to $\min\limits_{a} F(a)$(=0), while the minimizer or the parameter
estimator $a_E$  provides an accurate estimation of
the system parameter $a$. For sufficiently small $\varepsilon>0$,
the objective function value $\mathbb{F}(a)$ for the slow system
will also get
  closer and closer to $0$, and the minimizer or the parameter  estimator $a_E^S$ is a good approximation of $a_E$.

\end{example}

\newpage


\begin{thebibliography}{99}
\bibitem{Barati} R. Barati, Parameter Estimation of Nonlinear
Muskingum Models Using Nelder-Mead Simplex Algorithm, \emph{J.
Hydrol. Eng.}, 2011, \textbf{16},   946-954.

\bibitem{Barton}R. R. Barton and J. S. Ivey, Modifications of the
Nelder-Mead Simplex Method for Stochastic Modification Response
Optimization, \emph{Proceeding of the 1991 Winter Simulation
Conference}, 1991,   945-953.



\bibitem{Bishwal} J. P. N. Bishwal,
\newblock{Parameter Estimation in Stochastic Differential Equations},
Springer, New York, 2007.

\bibitem{Chang} K. Chang, Stochastic Nelder-Mead simplex method - A
new globally convergent direct search method for simulation
optimization, \emph{European Journal of Operational Research},
\textbf{220}(2012),   684-694.


\bibitem{FuLiuDuan} H. Fu, X. Liu and J. Duan,
Slow manifolds for multi-time-scale stochastic evolutionary systems.
 \emph{Comm. Math. Sci.}, 2013,  Vol. \textbf{11}, No. 1,   141-162.
arXiv:1205.0333 [math.PR]

\bibitem{Ibragimov} I. A. Ibragimov and R. Z. Has'minskii,
\newblock{\em Statistical Estimation---Asymptotic Theory}.
Springer, New York, 1981.


\bibitem{Kan2013} X. Kan, J. Duan,
I. G. Kevrekidis and A. J. Roberts, Simulating Stochastic
Inertial Manifolds by a Backward-Forwand Approach. \emph{SIAM J. on Applied Dynamical Systems}, in press, Vol. 12, 2013.   arXiv:1206.4954 [math.DS].

\bibitem{Lagarias}
J. C. Lagarias, J. A. Reeds, M. H. Wright and P.
E. Wright, Convergence Properties of the Nelder-Mead Simplex Method
in Low Dimensions, \emph{SIAM J. Optim.}, 1998, \textbf{9}(1), pp
112-147.

\bibitem{Nocedal} J. Nocedal and S. J. Wright, \emph{Numerical
Optimization}, Springer Science+Business Media, LLC, New York, 2006.

\bibitem{Pham} N. Pham and B. M. Wilamowski, Improved Nelder Mead's
Simplex Method and Applications, \emph{Journal of Computing}, 2011,
\textbf{3}(3),   55-63.

\bibitem{Price} C. J. Price, I. D. Coope and D. Byatt, A converging Variant
of the Nelder-Mead Algorithm, \emph{Journal of Optimization Theorem
and Applications},  \textbf{113}(1), 2002,  pp. 5-19.

\bibitem{RenDuanJones} J. Ren, J. Duan and C. K. R. T.
Jones, Approximation of Random Slow Manifolds and Settling of
Inertial Particles under Uncertainty. Submitted to \emph{J. Dynamics \& Diff. Eqns.}, 2012.   arXiv:1212.4216

\bibitem{Schm} B. Schmalfuss and K. R. Schneider,   Invariant manifolds for
random dynamical systems with slow and fast variables. \emph{J. Dynamics \& Diff. Eqns.} \textbf{20} (2008), No. 1, 133--164.



\bibitem{SunDuan1}  X. Sun, J. Duan and X. Li, An impact of noise on invariant manifolds in nonlinear dynamical systems,
\emph{J. Math. Phys.} \textbf{51}, 042702 (2010).

\bibitem{SunDuan2} X. Sun, X. Kan and J. Duan, Approximation of invariant foliations for stochastic   dynamical systems.
\emph{Stochastics and Dynamics} \textbf{12}, No. 1, (2012), 1150011.

\bibitem{JYangDuan} J. Yang and J. Duan, Quantifying Model
Uncertainties in Complex Systems. \emph{Progress in Probability},
2011, Vol.65,  49-80.

\end{thebibliography}
\end{document}